\documentclass{article}

\usepackage{PRIMEarxiv}
\newcommand{\inlineequation}[2][]{\begin{equation}\label{#1}#2\end{equation}}

\usepackage[utf8]{inputenc} 
\usepackage[T1]{fontenc}    
\usepackage{hyperref}       
\usepackage{url}            
\usepackage{booktabs}       
\usepackage{amsfonts}       
\usepackage{nicefrac}       
\usepackage{microtype}      
\usepackage{lipsum}		
\usepackage{graphicx}
\usepackage{doi}
\usepackage{amsmath}
\usepackage{lipsum}
\usepackage{amsfonts}
\usepackage{graphicx}
\usepackage{epstopdf}
\usepackage{algorithm} 
\usepackage{algpseudocode} 
\ifpdf
  \DeclareGraphicsExtensions{.eps,.pdf,.png,.jpg}
\else
  \DeclareGraphicsExtensions{.eps}
\fi

\newcommand{\G}{\mathcal{G}}
\newcommand{\bZero}{\vec{\boldsymbol{0}}}
\newcommand{\bmu}{\boldsymbol{\mu}}
\newcommand{\bP}{\boldsymbol{{u}}}

\newcommand{\boldx}{\boldsymbol{x}}

\DeclareMathOperator*{\argmin}{arg\,min}

\title{Aging and mortality of persons with HIV: a novel Kalman Filtering and DMD framework}

\author{Alex Viguerie\\
	Department of Pure and Applied Sciences\\
	Universit\`a degli studi di Urbino Carlo Bo\\
	Urbino, PU (Italy) 61029 \\
	\texttt{alexander.viguerie@uniurb.it} \\
	\And
        Elisa Iacomini	\\
        Department of Environmental and Prevention Sciences\\
	Universit\`a degli studi di Ferrara\\
	Ferrara, FE (Italy) 44121 \\
	\texttt{elisa.iacomini@unife.it} \\
}

\begin{document}
\maketitle

\begin{abstract}
Due to the widespread availability of effective antiretroviral therapy (ART) regimens, average lifespans of persons with HIV (PWH) in the United States have increased significantly in recent decades. In turn, the demographic profile of PWH has shifted. Older persons comprise an ever-increasing percentage of PWH, with this percentage expected to further increase in the coming years. This has profound implications for HIV treatment and care, as significant resources are required not only to manage HIV itself, but associated age-related comorbidities and health conditions that occur in aging PWH. Effective management of these challenges in the coming years requires accurate modeling of the PWH age structure. In the present work, we introduce several novel mathematical approaches related to this problem. We present a workflow combining a PDE model for the PWH population age structure, into which publicly-available HIV surveillance data is assimilated using the Ensemble Kalman Inversion (EKI) algorithm. This procedure allows us to rigorously reconstruct the age-dependent mortality trends for PWH over the last several decades. To project future trends, we introduce and analyze a novel variant of the Dynamic Mode Decomposition (DMD), nonnegative DMD. We show that nonnegative DMD provides physically-consistent projections of mortality and HIV diagnosis while remaining purely data-driven, and not requiring additional assumptions. We then combine these elements to provide forecasts for future trends in PWDH mortality and demographic evolution in the coming years.
\end{abstract}

\keywords{HIV\and population-structured models \and inverse ensemble Kalman filter\and Dynamic Mode Decomposition}

\section{Introduction}

\par
\begin{figure}
    \centering
    \includegraphics[width=\textwidth]{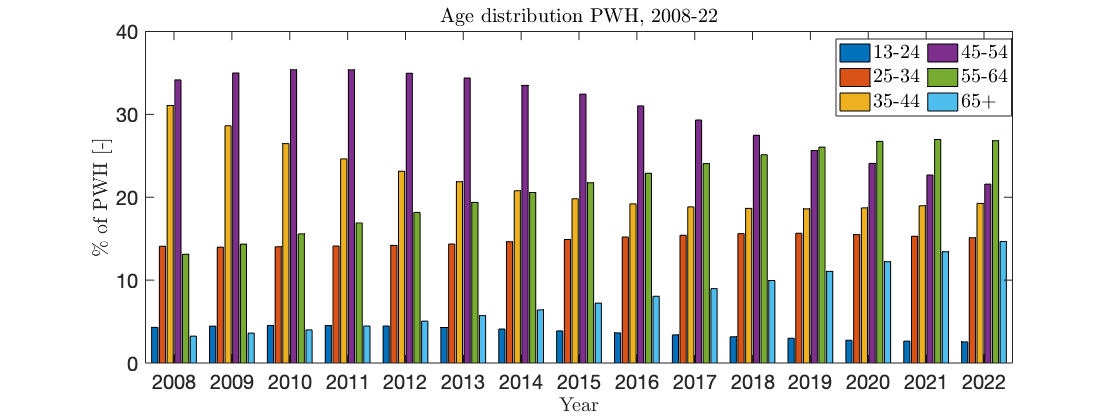}
    \caption{Age distribution of the diagnosed PWH population in the United States over the years 2008-22. } 
    \label{fig:PWHAgeStructure}
\end{figure}

 In recent decades, the life expectancy of individuals living with HIV, particularly in industrialized countries, has increased significantly, largely due to advancements in antiretroviral therapy (ART) \cite{trickey2023life, bosh2020vital}. ART has transformed HIV from a fatal disease into a manageable chronic condition, allowing many persons with HIV (PWH) to live lifespans nearly commensurate with the general population \cite{trickey2023life, marcus2020comparison,gueler2017life}. As a result, the demographic profile of persons with diagnosed HIV (PWDH) in the United States has shifted. As shown in Fig. \ref{fig:PWHAgeStructure} persons aged 55 and older have gone from approximately 16\% of the PWDH population in 2008 to nearly 45\% in 2022 \cite{ATLAS}. Such a drastic demographic shift presents new challenges in the fight against HIV, as long-term management of chronic HIV infection will require an increasing share of resources.

\par The growing population of older PWDH requires the healthcare system to address not only the ongoing management of HIV itself, but, increasingly, comorbidities and age-related conditions that arise as this population ages. These include cardiovascular disease, diabetes, osteoporosis, and certain cancers, which occur at higher levels, and at younger ages, among persons with HIV as compared to the general population \cite{bloch2020managing,bigna2020global,rodes2022ageing, roomaney2022aging, fosnacht2013older, collins2020means, nanditha2021excess,petoumenos2011hiv}. Older persons with HIV also may suffer from cognitive and functional impairment \cite{paul2018cognitive,mateen2012aging}. Additionally, the high burden of health comorbidities leads to high levels of polypharmacy among older PWDH \cite{smith2022polypharmacy}. Finally, although ART has significantly increased overall lifespan among PWDH, ART itself is associated with health complications arising from long-term use \cite{guaraldi2014morbidity,smith2022polypharmacy,chawla2018review}. As a result, healthcare providers must adapt their approaches and strategies to encompass a broader spectrum of health issues among persons with HIV \cite{collins2020means,kiplagat2022health}.
\par Given the changing landscape of HIV care, the focus of public health strategies must also evolve. Prevention efforts now must include the mitigation of health complications associated with HIV in older adults, long-term HIV infection, and long-term ART use among PWDH  \cite{collins2020means,kiplagat2022health}, necessitating a comprehensive approach that integrates HIV care with general healthcare for chronic disease prevention and management \cite{petoumenos2011hiv,kiplagat2022health,collins2020means}. Effective strategies will be crucial to maintain the health and quality of life of the aging PWDH, ensuring that the progress made in extending life expectancy is not undermined by the rise of preventable chronic conditions, and that increased life-years are not accompanied by decreased life quality. Furthermore, despite decreases in new HIV infections and diagnoses in recent years, improved life expectancy among PWDH has led to a significant increase in the PDWH population as a whole in the United States and other high-income countries \cite{cdcSurvReport2024, Puglia2024Toscana}, compounding aging-related difficulties, as HIV care must be allocated across a large, and increasing, population.
\par To help address these challenges, the development of accurate mathematical models capable of providing reliable projections for the PWDH age structure in coming years is an important priority in HIV policy research. Such models can play a fundamental role in the planning and evaluation of intervention strategies and the allocation of valuable HIV treatment and prevention services.  Recent work on this topic includes \cite{stover2024methods}, which applied a compartment model for demographic projections, and \cite{althoff2024forecasted, hyle2023growing}, which used an agent-based model to project how demographic changes in the US PWDH population may affect the burden of comorbidities. A particularly important component of these models is their ability to accurately characterize age-dependent PWDH mortality rates, and furthermore, how they may change in coming years. This arguably forms the foundation of such models, as the overall cost of HIV-related care depends most heavily on total PWDH person-years. 
\par Given the importance of the topic, in the present work, we seek to estimate how the age-dependent mortality rates for PWDH in the US have changed in recent years in a mathematically rigorous manner, and furthermore, extrapolate from these data to project future such changes. To the authors' knowledge, no such study has been undertaken, certainly not to describe changes in mortality at the age- and time-continuous level. To accomplish this task, we employ a data assimilation approach, in which age-discrete HIV surveillance data is integrated into a partial differential equation (PDE) model via an inverse ensemble Kalman filter (EnKF) to reconstruct the time evolution of age-dependent PWDH mortality since 2009. To project future trends in age-dependent mortality based on our reconstructions, we develop and apply a novel variant of Dynamic Mode Decomposition (DMD), nonnegative DMD. By incorporating these data into our PDE model, we provide a forecast PWDH age structure, annual HIV diagnoses, and annual PWDH deaths in the United States over the coming years.
\par To reconstruct the evolution of age-dependent PWDH mortality rates in recent years, we employ the Ensemble Kalman Inversion (EKI). The EKI is an adaptation of the ensemble Kalman filter (EnKF), a powerful tool for estimating the state of a dynamical system by iteratively refining predictions based on noisy observations. Originally introduced in the 1990s \cite{Evensen1994}, the EnKF has gained widespread use due to its ability to fuse model predictions with measurement data. Ensemble Kalman filtering has been employed to  enhance the accuracy of solutions in diverse applications, including oceanography \cite{EvensenVanLeeuwen1996}, reservoir modeling \cite{Aanonsen2009}, weather forecasting \cite{McLaughlin2014}, milling process \cite{SchwenzerViscontiEtAl2019}, geophysical applications \cite{Keller2021}, physics \cite{Li2020} and machine learning \cite{KovachkiStuart2019,Alper}. In public health, the EnKF was recently used to forecast heroin overdose deaths \cite{bottcher2024forecasting}. The EKI is a modified version of the EnKF designed for inverse problems, wherein underlying model parameters are estimated given a time-series of observations corresponding to model states \cite{schillings2017analysis,herty2019kinetic,herty13recent}. The EKI is attractive computationally, as it does not require the computation of derivatives of the underlying model, as well as mathematically, as it allows for provable convergence results \cite{schillings2017analysis,herty2019kinetic}. 
\par To project future trends from the reconstructed mortality data, we introduce a novel modification of the Dynamic Mode Decomposition (DMD), nonnegative DMD. DMD is a Scientific Machine Learning (SciML) technique that extracts the most relevant dynamical structures existent in time-series data using a purely data-driven approach, with applications ranging from short-time future estimates to control, modal analysis and dimension reduction \cite{Kutz2016book}. DMD has been deployed in a wide range of scientific and engineering applications such as biomechanics \cite{viguerie2022data}, oncology \cite{viguerie2022data} epidemiology \cite{Proctor2015, BGVRC2021, viguerie2022coupled}, climate modeling \cite{Kutz2016}, aeroelasticity \cite{Fonzi2020}, additive manufacturing \cite{viguerie2022coupled}, urban mobility \cite{Alla2020}, and the modeling and simulation of batteries \cite{alla2024piecewise}. The authors are not aware of any work using DMD to forecast the evolution of mortality rates, or similar quantities. In general, however,  standard DMD fail to respect non-negativity when used for forecasting beyond the training period, even when the underlying training data are nonnegative \cite{takeishi2017sparse, viguerie2022coupled}. This makes DMD in its original formulation potentially unsuitable for forecasting nonnegative quantities, such as mortality rates, and as such must be modified accordingly.


\par The article is outlined as follows. In the methods section, we provide a high-level description of the mathematical model and the EnKF method. We then introduce the proposed novel variant of DMD, nonnegative DMD, discuss its numerical solution, and establish some of its basic mathematical properties. We then present the results of our analyses. We first present the  reconstructions of age-dependent PWDH mortality curves over previous years. We then use nonnegative DMD to obtain projections of future mortality rates, and apply these rates in our underlying PDE model in order to provide projections for age-dependent PWDH prevalence and mortality through 2030. We then provide a detailed discussion of the implications of our findings and conclude by discussing possible future directions to extend this research.

\section{Methods}
\subsection{Mathematical model}
We consider an age-structured PDE model for the population $u(a,t)$, which gives the size of the PWDH population $u$ aged $a$ at a time $t$. The basic equations read:
\begin{align}\begin{split}
\dfrac{\partial u}{\partial t} + \dfrac{\partial u}{\partial a} &= -\mu(a,t)u + \lambda(a,t) \,\,\, \text{ for } t_0 < t\leq t_{end},\, 0\leq a < \infty,   \\
u(a,0) &= u_0(a) \,\, \text{ at } t=t_0, \,\,\, u(0,t)= g(t),
\end{split}\label{eq:model}\end{align}
where $\mu(a,t)\geq 0$ is the age- and time-dependent mortality rate, $u_0 \geq 0$ is the age-distribution of the PWDH population at the initial time, and $\lambda(a,t) \geq 0$ gives the number of new HIV diagnoses of persons aged $a$ at a time $t$, and is known from data. 
\par The inflow condition $g(t)\geq 0$ defines new births, and may be time- and/or state-dependent in general. In the present, however, we assume $g(t)=0$, for several reasons. First, levels of perinatal HIV in the US are very low \cite{lampe2023achieving}. Second, publicly available HIV surveillance data in the US does not include persons younger than 13 years of age \cite{cdcSurvReport2024}. Finally, the method of generating continuous age-distributions from discrete age-bins available from surveillance data (discussed in section 2.2), used to generate $u_0(a)$ and $\lambda(a,t)$, ensures $u(a,t)>0$ and remains smooth for all $a>0$ throughout the simulation period.

\par In practice, we solve \eqref{eq:model} numerically, considering a bounded age domain $0\leq a \leq 101$. We employ a second-order backward differentiation (BDF2) scheme throughout the age domain, assuming physically consistent `ghost node' values of $a_{-1}^t=a_{-2}^t =0$ for all $t$. For the time discretization, we use Heun's method to solve the first time step, and BDF2 for the remaining time steps. 
\newline Models of this type are common in many applications in the natural and social sciences, and we refer the reader to \cite{murray2003mathematical} for a presentation of these models and \cite{li2020age} for an extensive discussion of more complex extensions. 

\par Denote the domain of $u(\cdot,t)$, the set of admissible human ages, as $\mathbb{A} = \mathbb{R}^+ = \lbrack 0,\,\infty)$, endowed with a $\sigma$-algebra $\mathfrak{B}$ and nonnegative measure $\nu(A)$. In the following, we use $u(a,t)$ to denote a generic solution of \eqref{eq:model} and provide the following positivity result:
\newline \textbf{Theorem 1: }  \textit{ If $u_0\geq0$ for all $a$ and $\mu\geq 0,\, \lambda \geq 0 $ for all $(a,t)$, then $u(a,t)\geq0$ for all $(a,t)$.   }
\newline \textbf{\textit{Proof.}} Since the age $a$ advances in sync with time, we can consider the \textit{characteristic curve} $a(t)=a_0 + t$ for some fixed $a_0$, and note that $da/dt = 1$. Along the characteristic, we can rewrite $\eqref{eq:model}$ in terms of a function $v(t)$, depending only on $t$: 
\begin{equation}\label{eq:defCurveU}
v(t) := u(a(t), t) =u(a_0+t, t).
\end{equation}
Differentiating with respect to $t$ gives:
\begin{equation}\label{eq:curveDiff}
\dfrac{dv}{dt} = \dfrac{\partial u }{\partial t} + \dfrac{\partial u}{\partial a}\dfrac{d}{dt}\lbrack a_0 + t\rbrack = \dfrac{\partial u }{\partial t}(a_0+t,\,t) + \dfrac{\partial u}{\partial a}(a_0+t,\,t).
\end{equation}
where the last equation on the right hand side follows from $da/dt=1$. From \eqref{eq:defCurveU}, \eqref{eq:curveDiff} along the characteristic \eqref{eq:model} reduces to the ordinary differential equation:
\begin{equation}
    \dfrac{dv}{dt}= -\mu(a_0+t,t)v(t) + \lambda(a_0+t,t).
\end{equation}
Re-arranging, multiplying by the integrating factor $e^{\int_0^t\mu(a_0 +\tau,\tau)d\tau}$, and integrating from $0$ to $t$: 
\begin{equation}
 v(t)e^{\int_0^t\mu(a_0 +\tau,\tau)d\tau}  = v(0)+\int_0^{t} \lambda(a_0+\xi,\xi) e^{\int_0^\xi\mu(a_0 +\tau,\tau)d\tau}d\xi.
\end{equation}
Recalling that $v(0)=u_0(a_0)$:
\begin{equation}
u(a_0+t,t) =e^{-\int_0^t \mu(a_0+\tau,\tau) d\tau }\left(u_0(a_0) + \int_0^{t} \lambda(a_0+\xi,\xi) e^{\int_0^\xi\mu(a_0 +\tau,\tau)d\tau}d\xi\right).
\end{equation}
From our assumptions on $\mu,\,u_0,\,\lambda$, all terms on the right-hand side are nonnegative, and hence:
\begin{equation}\label{a0bound}
u(a_0+t,t)\geq 0 \,\,\forall t.
\end{equation}
Furthermore, since $a_0$ was arbitrary, \eqref{a0bound} holds for any $a_0 \geq 0 $, implying:
\begin{equation}\label{positivityResult}
u(a,t)\geq 0 \,\,\, \forall (a,t),
\end{equation}
as was to be shown.

\subsection{Ensemble Kalman Inversion}


The Ensemble Kalman Inversion (EKI) is a computational technique that extends the principles of the Kalman filter to solve inverse problems by utilizing an ensemble of model realizations. By iteratively updating this ensemble with observational data, EKI effectively estimates unknown parameters and reduces uncertainty, making it particularly suitable for high-dimensional and nonlinear inverse problems \cite{schillings2017analysis,herty2019kinetic,herty13recent}. 

Here, we describe the Ensemble Kalman Inversion procedure applied in the current work.
In order to do this, let us introduce the notation we will use in the following:
\begin{itemize}
    \item The overline symbol $\overline{\,\cdot\,}$ denotes an averaged quantity
 \item $\boldsymbol{u}(a,t)$ denotes the age-distribution of the PWDH population in time; 
 \item $\boldsymbol{\mu}(a,t)$ is our unknown; in this case, it defines the time- and age-dependent mortality curve;
 \item $\Gamma$ is the matrix of the covariance of the noise, which is assumed to be known;
    \item $\G(\cdot)$ is our output of the total annual PWDH deaths in each age bracket, calculated from the population state $\bP$ and mortality curve $\mu$.  
    \item $n_{kf}$ is the number of Kalman filtering steps.
    \item $A(\boldsymbol{\mu},\lambda)$ is the operator which advances population age distribution $\Delta t$ units in time for a given $\bmu(a,t)$ and  $\lambda(a,t)$ according to the system dynamics defined by \eqref{eq:model}:
    \begin{equation}
    \boldsymbol{u}(a,t+\Delta t) = A\left(\boldsymbol{\mu}(a,t),\lambda(a,t)\right)\boldsymbol{u}(a,t).
    \end{equation}
\end{itemize}
For the sake of readibility, we hereafter write $\bmu(a,t),\,\bP(a,t),\,$ and $\lambda(a,t)$ respectively as $\bmu_{t},\,\bP_{t},\,$ and $\lambda_{t}$, with the dependence of each quantity on age $a$ understood if not explicitly denoted.
\par We seek to use the EKI to estimate the time- and age-dependent mortality curve $\bmu_t$, which is not observable directly, by incorporating information from the model \eqref{eq:model} and surveillance data. The estimation of $\bmu_t$ relies on the current age distribution $P_t$, calculated from the \eqref{eq:model}, and the age-bracketed annual mortality data $y_t$, available from surveillance data \cite{ATLAS}. To ensure regularity in our estimate, we do not estimate $\bmu_t$ at each age individually, but compute its value at ages 1, 10, 22, 32, 42, 52, 57, 62, 67, 72, 77, 82, 95, 101. We then consider the full curve  $\bmu_t$ as the piecewise cubic Hermite interpolating polynomial through these points. We note that we performed a sensitivity analysis to the choice of interpolation points, and found that the effect on the overall reconstruction was not large, provided that there are sufficiently many points across all age ranges. 
\par We begin by initializing an ensemble of mortality curves at ${\bmu_{t,j}^1}$ for $j=1,\,...,\,J$ at $t=2009$, by sampling from the distribution of the mortality rate of the general (non-PWH) population $\bmu_0$. The superscript number refers to the current filter step, further clarified in the following. We assume a multivariate normal distribution with mean $\bmu_0$ and covariance $\text{diag}(\bmu_0^2)$, the component-wise square. 
\par Similarly, we also initialize an ensemble of population age structures $P_{t_0,j}$ for $j=1,\,2,\,...,\,J$ at the initial time $t_0=2009$, which serve as the initial conditions for \eqref{eq:model}. Each $P_{t_0,j}$ is generated from surveillance data for the diagnosed PWDH population at year-end 2008 \cite{ATLAS}. These data are stratified into age brackets: 13-24, 25-34, 35-44, 45-54, 55-59, 60-64, 65-69, 70-74, 75-79, 80-84, and 85+. We define each $P_{t,j}^1$ such that the overall number of PWDH in each age bracket is equal to the surveillance data, while the distribution of PWDH \textit{within} each age bracket are defined with a uniform sample. For each calendar year $t_i=2009-22,\,i=1,\,2,\,...$ during the simulation, we follow an identical procedure using the age-bracketed annual new HIV diagnosis data from \cite{ATLAS} to define an ensemble $\lambda_{t_i,j},\,j=1,\,2,\,...,\,J$ for new entries in the model \eqref{eq:model}.

\par Assume that \eqref{eq:model} is discretized time with time step $\Delta t= 1/m,\,m\in\mathbb{N}$, such that each year is divided into $m$ time-steps. For each simulated calendar year $t_i$, we perform $n_{kf}$ inverse Kalman Filtering steps. The procedure is given in Algorithm \eqref{alg:InvKF}. In a given year, we perform $n_{kf}$ filtering steps. At the $nn$-th filtering step, we first advance the $P_j^{nn}$ following \eqref{eq:model} by solving the model over $m$ time-steps, using $\lambda_{t_i,j}$ and the mortality curve $\bmu_{t_i,j}^{nn}$. Note that we must keep a running update of the simulated age-bracketed mortality $\G(P_j^{nn},\bmu_{t_i,j}^{nn})$ at each time step, so that an end-of-year total is available. 
\par If $nn$ is less than the total number of filter steps $n_{kf}$, we use the calculated $\G(P_j^{nn},\bmu_{t_i,j}^{nn})$ and the surveillance mortality data $y_{t_i}$ to obtain the updated mortality curve $\bmu_{t_i,j}^{nn+1}$. We then simulate the model for the year $t_i$ again, using with the updated $\bmu_{t_i,j}^{nn+1}$. If $nn=n_{kf}$, we advance to year $t_{i+1}$, setting $\bmu_{t_{i+1},j}^1 = \bmu_{t_i,j}^{n_{kf}}$ and generating $\lambda_{t_{i+1},j}$ from the surveillance data for year $t_{i+1}$. 
Note that the solution of the EnKF procedure is given by the mean of the ensemble at the step $n_{kf}$ for each year.
\newline
\newline \noindent \textbf{Remark.} Note that, as \eqref{eq:model} must be solved $n_{kf}$ times in each year, we must keep values of $P_j$ from the \textit{preceding} year in order to repeatedly solve the first time-step of the \textit{current} year. The exact number of time-step(s) to be saved depends on the specific time-stepping scheme used. As we used BDF2 in the current work, we kept the last two time-steps from the $n_{kf}$-th filtering step of the previous year.

\begin{algorithm}
	\caption{Algorithm with fixed number of filter steps for the model $\G$} \label{alg_gen}
	\begin{algorithmic}[1]
 \State Initialize state ensemble $\bP_{t_0,j}$
 \State Initialize measurement noise: $\boldsymbol{Q}=\text{diag}(\text{Var}(\G(\bP_{t_0},\bmu_0))$, a diagonal matrix consisting of the variance in observed deaths in each age group for the general-population mortality $\bmu_0$ over the initial state ensemble $\bP_{t_0}$.
	\State Initialize the mortality curve ensemble $\bmu_{t_0,j}^{n_{kf}} \sim N(\bmu_0,\,\text{diag}(\bmu_0^2))  $, with $j=1,\,2,\,...\dots,J$, where $\bmu_0^2$ is interpreted as the element-wise square, and $n_{kf}>0$ the number of optimization steps 
	\State Set $n=0$, $t_1=2009$ 
    \For{$t_i$=2009:2022}
        \State Generate ensemble of annual new diagnoses $\lambda_{t_i,j}$;
        \State Set $\bmu_{t_i,j}^{1}=\bmu_{t_{i-1},j}^{n_{kf}}$;
        \For{$nn$=1:$n_{kf}$} 
        \State Set $\G(\bP_j,\bmu_{t_i,j}^{nn})=\boldsymbol{0}$;
            \For{$k=1:m$}
            \State Advance the population structure ensemble and age-structured mortality output:
            \begin{align*}
                 \bP_{t_i + k\Delta t,\,j} &= A(\bmu_{t_i,j}^{nn}, \lambda_{t_i,j})\bP_{t_i+(k-1)\Delta t,\,j  }; \\
                 \G(\bP_j,\bmu_{t_i,j}^{nn}) &= \G(\bP_j,\bmu_{t_i,j}^{nn})+ \Delta t \,\text{diag}(\bmu_{t_i,j}^{nn})\,  \bP_{t_i + k\Delta t,\,j}; 
            \end{align*}
            \EndFor
			\State $\bar{\boldsymbol{\mu}}^{nn}=\frac{1}{J} \sum_{j=1}^J \bmu_{t_i,j}^{nn}$ 
			\State $\overline{\G}^{nn}=\frac{1}{J} \sum_{j=1}^J \G(\bP_j,\bmu_{t_i,j}^{nn})$ 
            \State $\boldsymbol{\eta} \sim \mathcal{N}(0,\boldsymbol{Q})$
   \State Solve the EnKF procedure:
			\begin{align*}	
   C(\boldsymbol{\mu}^{nn})&= \frac{1}{J}\sum_{j=1}^J(\bmu_j^{nn}-\overline{\boldsymbol{\mu}}^{nn})\otimes(\G(\bP_j,\bmu_j^{nn})-\overline{\G}^{nn})\\
   D(\boldsymbol{\mu}^{nn})&= \frac{1}{J}\sum_{j=1}^J(\G(\bP_j,\bmu_j^{nn})-\overline{\G}^{nn})\otimes(\G(\bP_j,\bmu_j^{nn})-\overline{\G}^{nn})\tag{EnKF}\\
\bmu_j^{nn+1}&=\bmu_j^{nn}+  C(\boldsymbol{\mu}^{nn}) \left( D(\boldsymbol{\mu}^{nn})+ \Gamma^{-1} \right)^{-1} \left( y_{t_i} + \boldsymbol{\eta} -\mathcal{G}(\bP_j,{\bmu}_j^{nn}) \right)
\end{align*}
   
		\EndFor
\EndFor
	\end{algorithmic}\label{alg:InvKF} 
\end{algorithm}

\subsection{Nonnegative DMD}
\par In the present work, we introduce a novel variant of DMD, nonnegative DMD, and use it to extrapolate the evolution of the age-dependent mortality rates (calculated using the Ensemble Kalman Inversion) and the age structure of new HIV diagnoses (estimated from surveillance data) through the year 2030. In the following, we will provide a brief introduction to the standard DMD algorithm, then we will introduce the nonnegative DMD algorithm, and provide some information regarding its calculation. 
\subsubsection{Standard DMD}
In this subsection, we provide a brief description of DMD. For the purposes of simplicity and exposition, we restrict our discussion to the original formulation of DMD as proposed in \cite{schmid2010dynamic}. We remark that this is a vast topic of research, and refer the reader to \cite{Kutz2016book} for a more comprehensive treatment of DMD. Suppose we have time series consisting of $n$ snapshots:
\begin{equation}\label{timeSeries}
\big\{\boldx_i\big\}_{i=1}^n,\qquad\,\boldx_i \in \mathbb{R}^m\,\,\,\,\forall i.\end{equation} 
We arrange these snapshots into two $m \times n-1$ matrices $X_1,\,X_2$ column-wise as follows:
\begin{equation}\label{xTimeSeries}
X_1=\begin{bmatrix} | & | & \, & | \\ \boldx_1 & \boldx_2 & ... &\boldx_{n-1} \\  | & | & \, & | \end{bmatrix},\qquad X_2=\begin{bmatrix} | & | & \, & | \\ \boldx_2 & \boldx_3 & ... &\boldx_{n} \\  | & | & \, & | \end{bmatrix}.
\end{equation}
DMD seeks to reconstruct the operator $A$ mapping $X_1$ to $X_2$, that is:
\begin{equation}\label{shiftOperator}
    A X_1 = X_2,
\end{equation}
which can be obtained as \cite{schmid2010dynamic}:
\begin{equation}\label{dmdStandard}
    A = X_2 \left(X_1\right)^{\dag},
\end{equation}
where $X_1^{\dag}$ is the Moore-Penrose pseudoinverse of $X_1$ (see e.g. \cite{golub2013matrix}). This formulation of $A$ solves the minimization problem:
\begin{equation}
\argmin_{A} \| A X_1 - X_2 \|_{F},
\end{equation}
with $\|\,\cdot\,\|_F$ denoting the Frobenius norm \cite{kutz2016dynamic}. Note that, in practice, $A$ can be large and dense, and is not typically computed directly as \eqref{dmdStandard}. Instead, algorithms which efficiently compute lower-rank approximations of $A$ can be leveraged to avoid ever forming $A$ explicitly (see e.g. \cite{Kutz2016book, schmid2010dynamic, viguerie2022data}).
\subsubsection{Nonnegative DMD algorithm}\label{nonNegDMDSection}
\par Being a completely data-driven algorithm, DMD, particularly when extrapolated past the training interval, will not necessarily respect the physics of the underlying systems. For example, in \cite{alla2024piecewise}, it was demonstrated that standard DMD can fail to predict periodic limit cycles. Similarly, in \cite{viguerie2022coupled}, it was shown that DMD, when formulated as \eqref{dmdStandard} (or with a lower-rank approximation), will conserve the mass of a mass-conservative system (assuming the data is also mass-conservative), even projected out in time; however, other key physical properties, in particular nonnegativity, were not preserved in the extrapolations.
\par In the present work, we employ DMD to project two strictly nonnegative physical quantities: future HIV diagnoses and future age-dependent mortality rates among PWH. While standard DMD reliably reproduces the training data in both instances, attempting to extrapolate beyond the training interval results in negative, oscillatory reconstructions, even over short projection windows (1-2 years).

\par In order to guarantee $(A^{n} \boldx)_i \geq 0,\,\,\forall \,i,\,n$, assuming a nonnegative $\boldx$, it is enough to guarantee that the matrix $A$ itself is nonnegative. The most obvious solution would be to find $A$ as in \eqref{dmdStandard}, and then use nonnegative matrix factorization to replace $A$ with an approximate nonnegative version $A\approx\widehat{A}=MN$, that is:

\begin{equation}\label{firstAttempt}
\argmin_{\widehat{A}=MN} \| MN - X_2(X_1)^{\dag} \|_F,\,\,\,M_{i,j}, N_{i,j} \geq 0\,\,\forall i,\,j.
\end{equation}
However, our attempts in this direction were unsuccessful, and approximations of $A$ obtained in this way did not provide useful forecasts.
\par As $\{ \boldx \}_{i=1}^n >0$ element-wise for all $i$, we also attempted to define:
\begin{equation}\label{secondAttempt}\widetilde{X}_1=\log(X_1),\,\,\widetilde{X}_2=\log(X_2),\, \widetilde{A}=\widetilde{X}_2(\widetilde{X}_1)^{\dag},\end{equation}
and proceeding as:
\begin{equation}\label{secondAttemptPart2}
    \boldx_{t+1} = \exp\left(\widetilde{A}\log\left(\boldx_{t}\right)\right).
\end{equation}
Note the logarithm and exponential functions in \eqref{secondAttempt}, \eqref{secondAttemptPart2} refer to the elementwise operations. While the $\boldx_i$ obtained in this manner remained nonnegative, and accurately reconstructed the training period, extrapolating beyond the training period resulted in oscillatory, nonphysical solutions, even for just a single time-step, making this approach unsuitable for forecasting, at least for the problems studied herein.

Instead, we formulate nonnegative DMD as the following:
\par\noindent\textbf{Definition.} \textit{Let $X\in \mathbb{R}^{m\times n}$ and let $X_1$, $X_2$ be snapshot matrices as defined in \eqref{xTimeSeries}. We define the nonnegative Dynamic Mode Decomposition as the operator $A$ satisfying:}
\begin{equation}\label{nonNegDMD}
\argmin_A \| AX_1 - X_2\|_F,\,\,\,A_{i,j}\geq 0\,\,\forall i,\,j.
\end{equation}
In general, the solutions to \eqref{firstAttempt} and \eqref{nonNegDMD} will not coincide. As the definition of the DMD operator \eqref{dmdStandard} is motivated by dynamics \eqref{shiftOperator}, we believe that \eqref{nonNegDMD}, defined directly in terms of the underlying time dynamics in the data, is the correct formulation for the applications in the present. 

\par The computation of \eqref{nonNegDMD} is nontrivial. Note that well-known algorithms exist to solve the related matrix-vector \textit{nonnegative least squares} (NNLS) problem:
\cite{lawson1995solving}:
\begin{equation}
\argmin_{\boldx} \| A\boldx - \boldsymbol{b} \|_2  \,\,\,x_{j}\geq 0\,\,\forall j.\end{equation}
Problems related to \eqref{nonNegDMD}, in which one searches for a nonnegative matrix, have been studied in \cite{van2004fast, nadisic2022matrix, kim2007sparse, kiSrDh12}. We remark that, in each of these cases, it is assumed that one has access to the matrix $A$ and is attempting to find $X_1$, that is solving:
\begin{equation}\label{nonNegTypical}
\argmin_{X_1} \| AX_1 - X_2\|_F,\,\,\,X_{1_{i,j}}\geq 0\,\,\forall i,\,j.
\end{equation}
The DMD problem is the opposite, since $X_1$ is known from the data and $A$ is our unknown. However, recalling that $\|M\|_F = \|M^T \|_F$ \cite{golub2013matrix} for any matrix $M$, 
the transposed problem:
\begin{equation}\label{nonNegDMDTranspose}
\argmin_{A^T} \|X_1^T A^T - X_2^T\|_F,\,\,\,A_{i,j}^T\geq 0\,\,\forall i,\,j
\end{equation}
is equivalent to \eqref{nonNegTypical}.
Note that:
\begin{align}\begin{split}\label{frobNorm}
\|X_1^T A^T - X_2^T \|_F^2 &= \sum_{k=1}^{m} \|X_1^T (A_{k,1:m})^T-(X_{2\,\,k,1:(n-1)})^T\|_2^2 
\end{split}\end{align}
Let $\text{Vec}(M)$ denote the vector obtained of stacking the successive columns of a matrix $M$ on top of one another. Then we can write \eqref{frobNorm} as:
\begin{align}\begin{split}\label{matVecFormln}
\|X_1^T A^T - X_2^T \|_F^2 &= \|(I \otimes X_1^T )\text{Vec}(A^T)  - \text{Vec}(X_2^T)\|_2^2,
\end{split}\end{align}
where $I$ is the identity matrix of appropriate dimension, $\otimes$ denotes the Kronecker product. Note that \eqref{matVecFormln} is a matrix-vector system. The problem \eqref{nonNegDMD} is therefore equivalent to the following matrix-vector NNLS problem:
\begin{equation}\label{nonNegDMDMatVec}
\argmin_{\text{Vec}(A^T)} \|(I \otimes X_1^T )\text{Vec}(A^T)  - \text{Vec}(X_2^T)\|_2,\,\,\,\text{Vec}(A^T)_j\geq 0\,\,\forall j.
\end{equation}

Standard NNLS algorithms (see e.g.  \cite{lawson1995solving}) can be applied on \eqref{nonNegDMDMatVec}. Converting $\text{Vec}(A^T)$ back into matrix form and transposing, one obtains a solution to \eqref{nonNegDMD}. For small problems, this approach is feasible, and it was used throughout the present. 
\par However, for high-dimensional data, or for datasets with many observations, the direct vectorization approach may be untenable. In \cite{nadisic2022matrix}, the authors discuss several methods for solving matrix-matrix NNLS problems, including techniques that are strictly equivalent to \eqref{nonNegDMDMatVec}, as well as techniques that seek sparse solutions to \eqref{nonNegTypical} directly. Note that the vectorized problem \eqref{nonNegDMDMatVec} has the structure:
\begin{equation}
\argmin_{A} \bigg\| \begin{pmatrix} X_1^T & 0 & ... & 0 \\ 
                        0 & X_1^T & ... & 0 \\
                        \vdots & \vdots & \ddots &\vdots \\
                         0 & 0 & ... & X_1^T \end{pmatrix} \begin{pmatrix}A_{(1,:)} \\ A_{(2,:)} \\ \vdots \\ A_{(m,:)}\end{pmatrix} -  \begin{pmatrix}X_{2\,(1,:)} \\ X_{2\,(2,:)} \\ \vdots \\ X_{2\,(m,:)}\end{pmatrix} \bigg\|_2^2, \end{equation}
and hence can be regarded $m$ independent NNLS problems:
\begin{equation}\label{separateColCalc}
    \argmin_{A_{(k,:)}} \| X_1^T A_{(k,:)} - X_{2\,(k,:)} \|_2^2,\,A_{(k,j)} \geq 0\,\,\forall \,j.
\end{equation}
Since $X_1,\,X_2$ are known, and each $A_{(k,:)}$ can be computed independently, \eqref{nonNegDMDMatVec} is easily parallelized. We remark that while the problems \eqref{nonNegDMDMatVec} and \eqref{separateColCalc} are equivalent in theory, in practice, solutions to NNLS for underdetermined systems are not unique in general, and must be obtained through iterative procedures. The solutions obtained when solving the formulations \eqref{nonNegDMDMatVec} and \eqref{separateColCalc} will therefore differ in practice. In the present, the authors observed that such differences were trivial; however, a comprehensive investigation of this issue was not pursued.
\par  We remark that related work on a nonnegative DMD variant was pursued in \cite{takeishi2017sparse}. The formulation introduced in the present differs from \cite{takeishi2017sparse} in that we require that the approximated \textit{operator} $A$ be strictly nonnegative, whereas the approach in \cite{takeishi2017sparse} requires only that the \textit{dynamic modes}, or DMD eigenvectors, are nonnegative. However, as nonnegativity is not required for the corresponding DMD eigenvalues, the reconstructed DMD operator is not nonnegative in general. We note that this difference in formulation reflects the primary intended usage of each algorithm, as the approach in \cite{takeishi2017sparse} was developed principally for diagnostic purposes and dimension reduction, whereas the present algorithm is intended for usage in forecasting and extrapolation.
\subsubsection{Application to PWDH age structure model}\label{sec:application}
We use the time- and age-dependent mortality curves obtained with the EKI for the years 2009-22 to forecast the evolution of the PWDH population over the years 2023-30. For the years 2009-22, we solve the model \eqref{eq:model}, directly applying the $\mu(a,t)$ obtained through the EKI procedure. 
\par To obtain mortality rates for the years 2023-30, we apply the nonnegative DMD algorithm outlined in section \ref{nonNegDMDSection} defined as:

\begin{align}\label{mortRateScen}
\boldsymbol{M} &= \argmin_{\boldsymbol{M}} \| \boldsymbol{M} \lbrack \bmu_{2009}\,\bmu_{2010}\,...\,\bmu_{2018} \rbrack -  \lbrack \bmu_{2010}\,\bmu_{2011}\,...\,\bmu_{2019} \rbrack \|_F,\,\, \boldsymbol{M} \geq 0, 
\end{align}
where the notation $\boldsymbol{M} \geq 0$ indicates that $\boldsymbol{M}$ is entry-wise nonnegative. The 2023 mortality rate is then defined as $\boldsymbol{M} \bmu_{2022}$, with the following years obtained by iterated multiplication by the $\boldsymbol{M}$. We did not include the years 2020-22 in training the nonnegative DMD operator, as the jump in mortality rate due to COVID-19 likely represents a temporary jump in mortality, rather than a long-term trend, as indicated by preliminary 2023 mortality data returning to pre-COVID levels \cite{ahmad2024mortality}.  
\par Similarly, we use nonnegative DMD  to project age-structured annual new HIV diagnoses over the years 2023-30 as:
\begin{equation}\label{diagRateScen}
    \boldsymbol{\Lambda} = \argmin_{\boldsymbol{\Lambda}} \| \boldsymbol{\Lambda} \lbrack \lambda_{2009}\,\lambda_{2010}\,...\,\lambda_{2018} \rbrack -  \lbrack \lambda_{2010}\,\lambda_{2011}\,...\,\lambda_{2019} \rbrack \|_F,\,\,\boldsymbol{\Lambda}\geq 0.
\end{equation}
As with mortality rates, we disregard the years 2020-22 in projecting HIV diagnoses, as available evidence suggests that the large drop in HIV diagnoses observed in 2020, followed by rebounds in 2021 and 2020, was primarily caused by disruption and subsequent recovery in HIV testing and diagnosis \cite{viguerie2023isolating,viguerie2024covid,dinenno2022hiv}. However, estimated underlying HIV incidence continued to follow the decreasing trends observed in the years immediately preceding the COVID-19 pandemic \cite{cdcSurvReport2024}. Therefore, the increases in HIV diagnoses observed in 2021 and 2022 likely reflect a recovery of diagnoses missed during 2020 due to reduced testing, rather than changes in HIV transmission. Hence, the pre-2019 trends in new HIV diagnoses are likely more reliable for informing future trends in new HIV diagnoses over the coming years.
\subsubsection{Mathematical analysis}
We provide a series of results that establish the key properties of nonnegative DMD. Although the results are presented in some generality, we assume some additional conditions on the underlying data, motivated by the specific problem studied in the present. Accordingly, they may not apply universally. 
\par Letting $X_1,\,X_2 \,\in \mathbb{R}^{m\times n-1}$ be as in \eqref{xTimeSeries}, we assume additionally:
\begin{itemize}
\item \textit{Strict positivity: }\inlineequation[strictPos]{(X_1)_{i,j}>0,\,(X_2)_{i,j}>0\,\,\forall \,i,j}
\item \textit{Temporal regularity: }\inlineequation[tempReg]{(X_1)_{i,j} < 2 (X_2)_{i,j},\,\, (X_2)_{i,j} < 2(X_1)_{i,j}}.
\end{itemize}
Both assumptions hold for the data $\bmu_{2009-19},\,\lambda_{2009-19}$ used for \eqref{mortRateScen}, \eqref{diagRateScen}, respectively.  
The meaning of \eqref{strictPos} is straightforward, and a similar assumption was employed in Theorem 1 in Section 2.1.1. Assumption \eqref{tempReg} states that the local temporal variation at each observation point (age in the present document) must remain bounded, and should not more than double, or less than halve, from time-step to time-step (the factor of two is chosen for convenience). As DMD assumes underlying continuity in time \cite{schmid2010dynamic, kutz2016dynamic}, \eqref{tempReg} should hold for any data sufficiently well-resolved in time; in the absence of noise, it should be possible in principle, to obtain a time series satisfying \eqref{tempReg} by increasing time-resolution of the measurements. However, we recognize that, in practice, time-series data is subject to limitations and \eqref{tempReg} may not hold in general. 
\par\noindent\textbf{Theorem 2. } \textit{Any $\boldsymbol{M}$ solving \eqref{mortRateScen}, or $\boldsymbol{\Lambda}$ solving \eqref{diagRateScen}, is nonzero.}
\par\noindent\textbf{Proof. } We provide a result for Problem \eqref{mortRateScen}, and note that an identical argument applies to \eqref{diagRateScen}. We prove contrapositive of the theorem statement, and show that if $\boldsymbol{M}$ is a matrix containing only zeroes, it cannot be a solution of \eqref{mortRateScen}. For ease of notation, let: 
\begin{equation}
 \widehat{\bmu}_{2009-18} = \lbrack \bmu_{2009}\,\bmu_{2010}\,...\,\bmu_{2018} \rbrack,\, \qquad  \widehat{\bmu}_{2010-19}= \lbrack \bmu_{2010}\,\bmu_{2011}\,...\,\bmu_{2019} \rbrack,
\end{equation}
and proceed by contradiction. Denote as $\bZero \in \mathbb{R}^{m\times m}$ the $m \times m$ matrix whose entries are identically zero, and suppose that $\bZero$ solves \eqref{mortRateScen}. Then:
\begin{equation}\label{frobNormIneq}
 \| \bZero  \widehat{\bmu}_{2009-18} -  \widehat{\bmu}_{2010-19}   \|_F = \|    \widehat{\bmu}_{2010-19}  \|_F  
 \leq   \| \boldsymbol{A}  \widehat{\bmu}_{2009-18} -  \widehat{\bmu}_{2010-19}   \|_F \, \forall \boldsymbol{A}\geq 0.  \end{equation}
Let $\boldsymbol{I}$ be the $m\times m$ identity matrix. Clearly $\boldsymbol{I}\geq 0$. Observe that: 
\begin{equation}\label{identityExpression}
 \| \boldsymbol{I}   \widehat{\bmu}_{2009-18} -   \widehat{\bmu}_{2010-19} \|_F = \|  \widehat{\bmu}_{2009-18} -   \widehat{\bmu}_{2010-19} \|_F.
\end{equation}
Squaring the right-hand side of \eqref{identityExpression} and applying the definition of the Frobenius norm gives: 

\begin{align}\begin{split}\label{manipulation}
   \sum_{j=2009}^{2018} \| \bmu_{j} - \bmu_{j+1} \|_2^2 
 &= \sum_{j=2009}^{2018} \sum_{i=1}^m (\mu_{i,j} - \mu_{i,j+1})^2 \\
 &= \sum_{j=2009}^{2018} \sum_{i=1}^m \mu_{i,j}^2 + \mu_{i,j+1}^2 - 2\mu_{i,j}\mu_{i,j+1} \\
 &< \sum_{j=2009}^{2018} \sum_{i=1}^m 2 \mu_{i,j}\mu_{i,j+1} + \mu_{i,j+1}^2 - 2\mu_{i,j}\mu_{i,j+1} \\ 
&=  \sum_{j=2009}^{2018} \sum_{i=1}^m \mu_{i,j+1}^2 \\ 
&= \|  \widehat{\bmu}_{2010-19} \|_F^2,
\end{split}\end{align}
where the third line follows from the temporal regularity assumption \eqref{tempReg}. From \eqref{manipulation}, it follows that:
\begin{equation}
    \| \boldsymbol{I} \bmu_{2009-18} - \bmu_{2010-19} \|_F < \| \bmu_{2010-19} \|_F =  \| \bZero\, \bmu_{2009-18} - \bmu_{2010-19} \|_F,
\end{equation}
contradicting \eqref{frobNormIneq}, implying that the zero matrix cannot be a solution of \eqref{diagRateScen}, as was to be shown.

\par\noindent\textbf{Interpretation in terms of the Perron-Frobenius operator.} In this portion of the article, we provide an interpretation of nonnegative DMD in terms of the Perron-Frobenius (P-F) operator. We note that numerical methods for the Perron-Frobenius operator have been studied in other settings \cite{klus2016numerical, klus2016Tensor, ding1993high, bollt2013applied}. DMD is generally interpreted as a numerical approximation of the \textit{Koopman operator}, the adjoint of the Perron-Frobenius operator \cite{kutz2016dynamic, brunton2022modern}. Roughly speaking, the Koopman operator acts on functions of the system state (commonly called \textit{observables}), while the P-F operator acts on \textit{densities} of the system state.
\par Since the P-F and Koopman operators are adjoint, one can use either to describe a given dynamical system \cite{lasota2013chaos}. Nonetheless, the literature on the numerical approximation of the two operators is somewhat distinct. Literature on the approximation of the P-F operator tends to focus on \textit{Ulam's method} and its variants \cite{klus2016numerical, ding1993high, bollt2013applied, ermann2010ulam, junge2009discretization, ding1991markov}, while approximation of the Koopman operator focuses on DMD and its variants \cite{brunton2022modern, kutz2016dynamic, alla2024piecewise, Proctor2015, colbrook2023residual}. Note that this statement is only a broad, general characterization, and several works discuss the approximation of both operators \cite{klus2016numerical, klus2020data, klus2018data} and/or hybrid-type approaches incorporating aspects of both DMD and Ulam's method \cite{goswami2018constrained, colbrook2023mpedmd}.

In the following, we describe why, in the context of the present work, and problem \eqref{mortRateScen} in particular, the P-F operator provides a natural interpretation. We remark that a full description of the Perron-Frobenius and Koopman operators is beyond the scope of the current work, and we recommend that the reader consult \cite{klus2020data} for a comprehensive discussion of the numerical approximation of the P-F and Koopman operators, and \cite{lasota2013chaos} for the underlying theory. We note that this description is not, nor is intended to be, fully rigorous. Rather, it is intended as an intuitive explanation for the nonnegative DMD formulation used herein, and how it differs from other common DMD approaches within the context of Koopman/P-F theory.
\par Consider a measure space defined by the three-tuple $\mathcal{M}:=(\mathbb{A},\mathfrak{B}, \nu)$ where \newline $\mathbb{A}=\mathbb{R}^+ = \lbrack 0,\,\infty)$, the nonnegative real numbers, is the space of possible ages, $\mathfrak{B}$ is a $\sigma-$algebra on $\mathbb{A}$ and $\nu$ is a $\sigma-$finite (positive) measure on $\mathbb{A}$, hereafter assumed to be the standard Borel measure unless otherwise specified. Let $S:\mathbb{A}\to\mathbb{A}$ be a measurable mapping which advances age of the population. For simplicity and without loss of generality, we assume $S$ advances age forward one unit; hence $S$ moves a member of the population aged $a$ to age $a+1$. For simplicity, we disregard new entries $\lambda(a,t)$ for the moment.
\par Let $u_t(a) \in \mathcal{L}^1(\mathbb{A},\,\nu)$ describe the population distribution at time $t$. We define the \textit{support} of $u_t$ as:
\begin{equation}
    \text{supp} (u_t):=  \{ a\in \mathbb{A} | u_t(a)\neq 0\}. 
\end{equation}
Since human lifespans are finite, we may assume:
\begin{equation}
\nu(    \text{supp} (u_t))<\infty.
\end{equation}
at all $t$. 
\par Since $u_t(a) \in \mathcal{L}^1(\mathbb{A},\,\nu),\,u_t(a)\geq 0 \,\forall \,a\in\mathbb{A}$, we may also define the \textit{population measure} $\widehat{u}_t\in\mathcal{M}$:
\begin{equation}\label{popMeasDef}
\widehat{u}_t(X) = \int_X u_t(a) da, \,\,\,X\in\mathfrak{B}.
\end{equation}
Define the three-tuple $\widehat{\mathcal{M}}:= (\mathbb{A},\,\mathfrak{B},\,\widehat{u}_t)$  and let  $\widehat{\mu}_t(X) \in \widehat{\mathcal{M}}$ be the \textit{mortality measure} at a time $t$, measuring the number of deaths in the age-set $X$ of the population $u$ at time $t$. Clearly, if $\widehat{u}_t(X)=0$ for some $X\in \mathfrak{B}$, then $\widehat{\mu}_t(X)=0$ and hence, by the Radon-Nikodym theorem \cite{lasota2013chaos}, there exists a density function $\mu_t \in \mathcal{L}^1(\mathbb{A},\,\widehat{u}_t)$ such that:
\begin{equation}\label{mortMeasDef}
\widehat{\mu}_t(X) = \int_X \mu_t(a) \widehat{u}_t(da) = \int_X \mu_t(a) u_t(a) da, \,\,\,X\in\mathfrak{B}.
\end{equation}

From the Riesz Representation Theorem \cite{lasota2013chaos}, there is a linear operator $\psi_t(u)$ on the dual space of $\mathcal{L}^1(\mathbb{A},\nu)$ such that for any $u\in\mathcal{L}^1(\mathbb{A},\nu)$:
\begin{equation}\label{RRT}
\psi_t(u) := \int_0^\infty u(a) \widehat{\mu}_t(da) = \int_0^\infty u(a) \mu_t(a) da = (u,\,\mu_t)= \widehat{\mu}_t(\mathbb{A}),
\end{equation}
where $(\cdot,\,\cdot)$ indicates the scalar product.
\par We now recall that our nonnegative DMD operator $A$ is recovered as the solution of:
\begin{equation}\label{nonNegDMDTranspose2}
\argmin_{A^T} \|X_1^T A^T - X_2^T\|_F,\,\,\,A_{i,j}^T\geq 0\,\,\forall i,\,j.
\end{equation}
For the sake of concreteness, suppose we are considering the mortality rate-recovery problem \eqref{mortRateScen}. Then the matrices $X_1$ and $X_2$ have the structures:
\begin{equation}\label{mortRatesDisc}
X_1^T = \begin{pmatrix} \bmu_{2009}^T \\ \bmu_{2010}^T \\ \vdots \\ \bmu_{2018}^T \end{pmatrix}, \qquad X_2^T = \begin{pmatrix} \bmu_{2010}^T \\ \bmu_{2011}^T \\ \vdots \\ \bmu_{2019}^T \end{pmatrix}.
\end{equation}
Intuitively, we may regard the discrete expression $\bmu_{t}^T u$ as an analogue of:
\begin{equation}\label{intuitiveExpl1}
    \bmu_t^T u \approx \int_0^\infty u(a) \widehat{\mu}_t(da) = (u,\,\mu_t)= \psi_t(u).
\end{equation}
Following a similar intuition as used in \eqref{intuitiveExpl1}, we note that the expression:
\begin{equation}
(A X_1 )^T u = X_2^Tu 
\end{equation}
is a discrete analogue of:
\begin{equation}
(u,\,A\mu_t) = (u,\,\mu_{t+1})
\end{equation}
which, together with \eqref{RRT}: 
\begin{equation}
\psi_{t+1}(u)= (u,\mu_{t+1}) = (u,A\mu_t) = (A^* u ,\mu_t) =\psi_{t}(A^* u ) = A\psi_t (u).
\end{equation}
Consequently, the nonnegative DMD operator can be interpreted as an approximation of the \textit{Perron-Frobenius operator} applied to the mortality measures $\widehat{\mu}_t$ (or their related densities concerning the population measure), thereby evolving them over time \cite{lasota2013chaos}. By employing NNLS to compute $A$, we ensure that the projected $\mu_t$ remain on the positive cone and in $\mathcal{L}^1(\mathbb{A},\widehat{u}_t)$ for every $t$. In other words, non-negative DMD guarantees that densities are transformed into densities, as desired.

\section{Results}
\subsection{Reconstruction of age-dependent PWDH mortality curves, 2009-22}
\begin{figure}
    \centering
    \includegraphics[width=\textwidth]{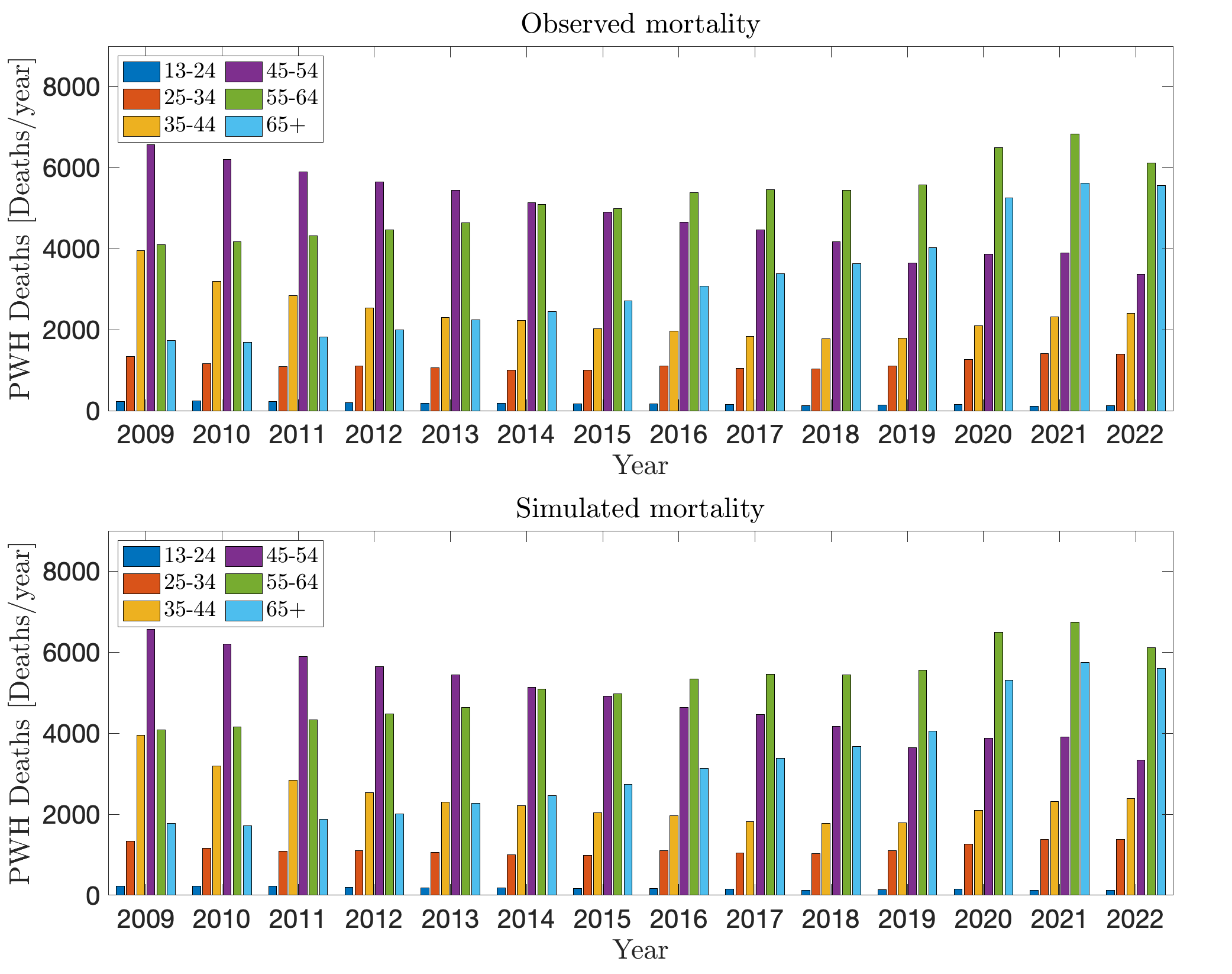}
    \caption{The comparison of surveillance (top) and simulated (bottom) mortality. }
    \label{fig:deathsBarChart}
\end{figure}

By applying the Kalman Inversion procedure outlined in the previous section, we observe a strong agreement with the data in Fig. \ref{fig:deathsBarChart}.

 This figure illustrates the age distribution of mortality for each year, juxtaposing real data (top panel) with simulated mortality (bottom panel). The simulated and observed mortality are in agreement both qualitatively and quantitatively.  
 \par
\begin{figure}
    \centering
    \includegraphics[width=\textwidth]{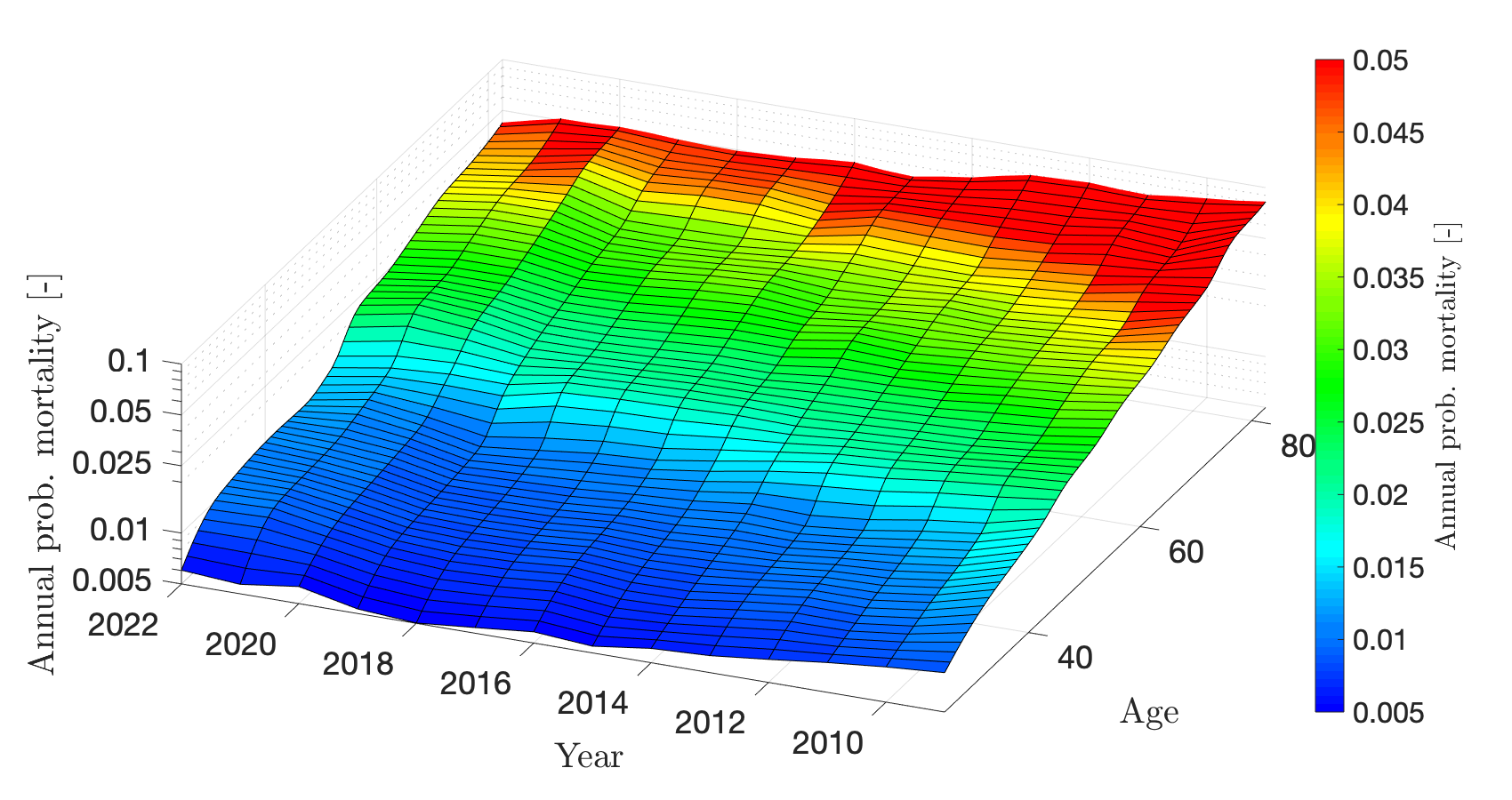}
    \caption{Annual probability of mortality at a given age, in time. Note the consistent decrease in time, punctuated by increases in the upper age ranges in 2020-22 caused by the COVID-19 pandemic. }
    \label{fig:deathSurface}
\end{figure}

 In Fig. \ref{fig:deathSurface}, we depict the annual mortality rate per age class over time. The color scale ranges from $0.005$ (dark blue) to $0.04$ (red). This visualization indicates a trend of decreasing annual mortality, evidenced by the expanding blue regions from 2010 to 2022. However, we note that at the beginning of the COVID-19 pandemic (2020), this trend slows somewhat, and mortality rates increases slightly. The red and yellow regions, which correspond to mortality rates at older ages, also show a strong diminishing trend in time until the beginning of the COVID-19 pandemic. The results suggest that the COVID-19 pandemic caused mortality rates across the 55+ age cohorts over the years 2020-21 to revert approximately to their approximate 2013 levels    \cite{ATLAS,viguerie2022mortality}.

\par
\begin{figure}
    \centering
    \includegraphics[width=\textwidth]{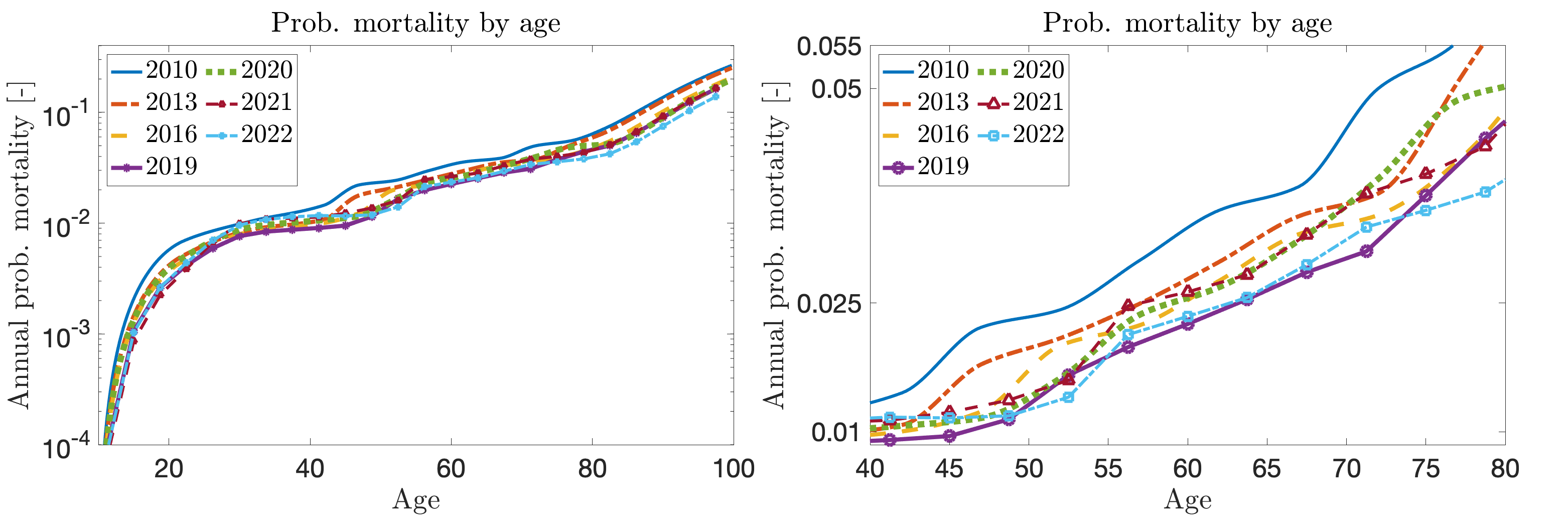}
    \caption{Mortality curves by age for several years, plotted side-by-side. The left panel shows mortality probability over the entire age range; the right panel focuses more closely on the important 40-80 age range. }
    \label{fig:mortCurves}
\end{figure}

 In Figure \ref{fig:mortCurves}, we depict several age-dependent mortality curves for individual years, over the entire age range (left panel), and focusing on the important 40-80 age range (right panel). Mortality decreases throughout the entire age range; however, as noted previously, these trends reverse in 2020, and mortality rates over the years 2020-21 increase among the 55+ cohort to levels similar to 2013-16. The 2022 rates show the effects of COVID subsiding among the 55+ cohort, with mortality rates consistent with levels observed during the years 2016-19 among PWDH aged 55-75, and at their lowest rates ever for PWDH aged 75+. However, we note that, even during the peak years of COVID-19 in 2020-21, mortality among older PWDH was still lower than 2010 levels, which suggests that the effects driving the decreases in mortality risk among PWDH over the past 15 years were more significant, on net, than the effects of the COVID-19 pandemic on PWDH mortality.

\subsection{Forecasting PWDH mortality, new diagnoses and PWDH population age-structure, 2023-30}
\subsubsection{Future PWDH mortality}
\begin{figure}
    \centering
    \includegraphics[width=\textwidth]{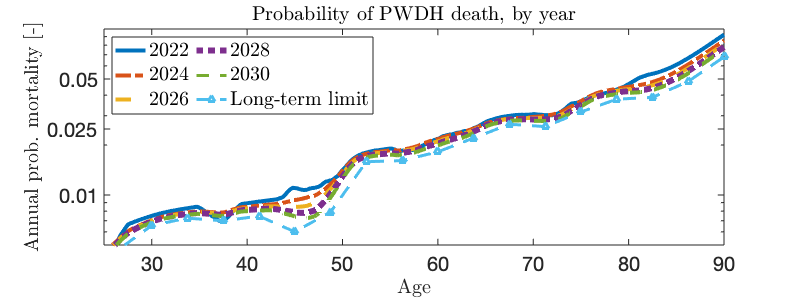}
    \caption{Projection of age-dependent mortality through 2030. The maximum eigenvector is the asymptotic limit of the projected future age-dependent PWDH mortality.}
    \label{fig:mortRateProject}
\end{figure}

In order to forecast the evolution of the mortality rate per age-class, we employ the nonnegative DMD as described in Sec. \ref{sec:application}. In Fig. \ref{fig:mortRateProject}, we depict the annual age-dependent probability of mortality for the years 2022, 2024, 2026, 2028, and 2030. The nonnegative DMD algorithm forecasts mortality rates to continue decreasing substantially in the 40-55 age range, with annual mortality projected to decrease 15.6\% in this age group from 2022 to 2030. Among PWDH aged 55-75, we also project future decreases in mortality, however these are comparatively modest, decreasing approximately 7.5\% from 2022 to 2030. For PWDH aged 75 to 90, mortality is projected to decrease approximately 15.6\% from 2022 to 2030. 
\par Analyzing the spectrum of $\boldsymbol{M}$, we find that it has a unique maximum eigenvalue $\xi_1 = \rho(\boldsymbol{M})= 1$ (consistent with the interpretation of $\boldsymbol{M}$ as an approximation of the Perron-Frobenius operator \cite{klus2020data, lasota2013chaos}), with second-largest eigenvalue $|\xi_2|=0.87$. We denote the corresponding eigenvectors as $\boldsymbol{v}_1$ and $\boldsymbol{v}_2$, respectively, and express $\bmu_{2022}$ as a linear combination of the eigenvectors $\boldsymbol{v}_i$ of $\boldsymbol{M}$:
\begin{equation}
    \bmu_{2022}=\sum_i a_i \boldsymbol{v}_i.
\end{equation}
Observe that:
\begin{equation}\label{mortLim}
    \lim_{k\to \infty } \boldsymbol{M}^k\bmu_{2022} = 1^k a_1 \boldsymbol{v}_1 + 0.87^k  a_2 \boldsymbol{v}_2 + \sum_{i\geq3} \xi_i^k a_i  \boldsymbol{v}_i = a_1 \boldsymbol{v}_1,
\end{equation}
since $|\xi_i| \leq .8909 < 1,\,\forall i\geq 3$. 
\par The eigenvector $a_1 \boldsymbol{v}_1$ of $\boldsymbol{M}$ corresponding to $\xi_1=1$ thus gives a long-term projection for the asymptotic age-dependent PWDH mortality rate, with age-dependent PWDH mortality converging to $a_1\boldsymbol{v}_1$ like $.87^t$. We plot this eigenvector as `Long-term limit' in Fig. \ref{fig:mortRateProject}. As the figure shows, beyond 2030, mortality among PWDH is projected to decrease significantly among PWDH aged 35-55 older than 80, and moderately among PWDH aged 60-80. We note that these forecasts are based on pre-2019 trends in PWDH mortality and hence cannot account for any changes beyond those pre-existing trends. As such, these forecasts, especially the asymptotic limit, should be interpreted with caution.
\par \textbf{Remark. } As $\boldsymbol{M}$ is nonnegative, if it is additionally irreducible, then we can guarantee that $\boldsymbol{v}_1$ has strictly positive entries from the Perron-Frobenius theorem. However, establishing irreducibility  for a general $\boldsymbol{M}$ is difficult, given the lack of a closed-form expression. Nonetheless, our experiments over an ensemble of $\boldsymbol{M}$ consistently showed that $\boldsymbol{M}$ has a unique maximum eigenvalue of 1, with a strictly positive corresponding eigenvector (we note that the spectral gap varied across the different $\boldsymbol{M}$). 
\subsubsection{Future HIV diagnoses}
\begin{figure}
    \centering
    \includegraphics[width=\textwidth]{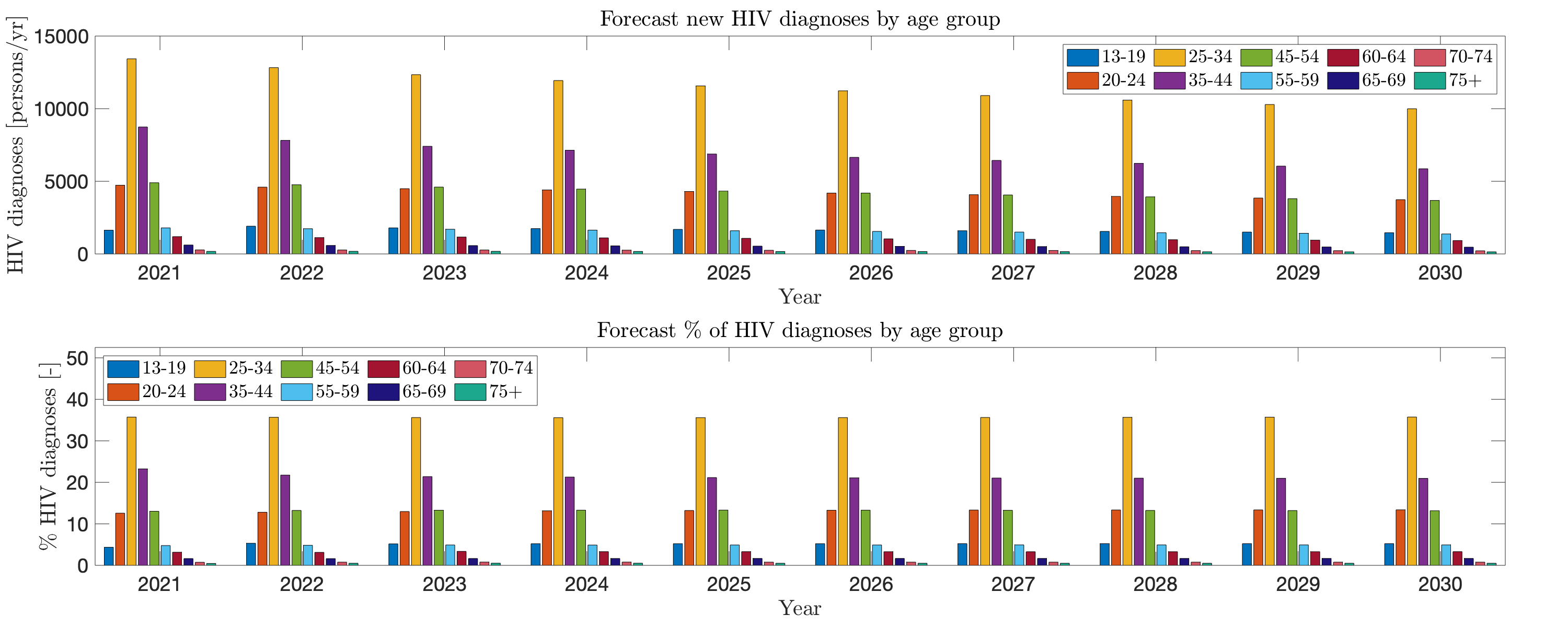}
    \caption{Projection of age-dependent HIV diagnoses through 2030. We project a slight decrease in overall diagnoses, with approximately 28,000 new HIV diagnoses in 2030 (compared to 34,500 in 2023). }
    \label{fig:diagProject}
\end{figure}

To forecast future new HIV diagnoses by age, we employed the nonnegative DMD algorithm on the age-depdendent diagnosis data as in \eqref{diagRateScen}. These data were defined using NCHHSTP surveillance data \cite{ATLAS}. New HIV diagnoses are projected to decrease at approximately 3\% per year, declining to approximately 28,000 new diagnoses in 2030. Furthermore, the age structure of new HIV diagnoses is not projected to show significant change in the coming years. We depict both the number and percentage of new HIV diagnoses by age group in Fig. \ref{fig:diagProject}.

\subsubsection{PWDH population age structure, 2023-30}
Using projected PWDH mortality \eqref{mortRateScen} and projected new HIV diagnoses \eqref{diagRateScen}, we solve the system \eqref{eq:model} through year-end 2030 to project the evolution of the PWDH age-structure over the years 2023-30. 
\par
\begin{figure}
    \centering
    \includegraphics[width=\textwidth]{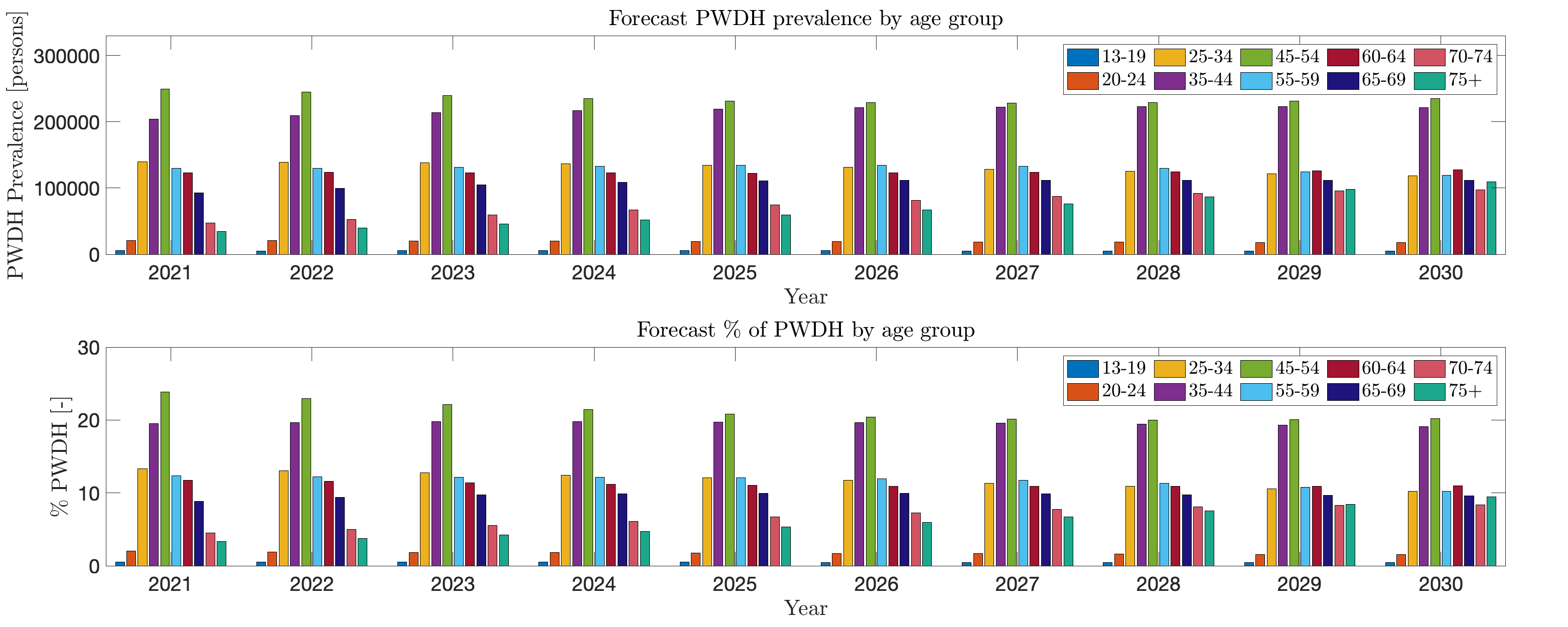}
    \caption{Projection of the PWDH age structure in the United States through 2030. We observe a gradual aging of the PWDH population. By 2030, 48.6\% of the PWDH population is projected to be over 55 and 27.4\% over 65.}
    \label{fig:PWHPrev1}
\end{figure}

 The total PWDH population is projected to grow over the simulation period, increasing to 1,160,000 by the end of 2030. As a whole, the PWDH population is also projected to age significantly and rapidly in the coming years. At year-end 2024, an estimated 42.5\% of PWDH were 55 and older, 19.4\% 65 and older, and 4.2\% 75 and older. By year-end 2030, we project the percentage of PWDH over 55 to increase to 47.\% and those over 65 to 26.1\%. The percentage of PWDH over 75 is projected to more than double over a five year span, reaching  8.5\%, by the end of 2030. This information is depicted in Fig. \ref{fig:PWHPrev1}.
\par
\begin{figure}
    \centering
    \includegraphics[width=\textwidth]{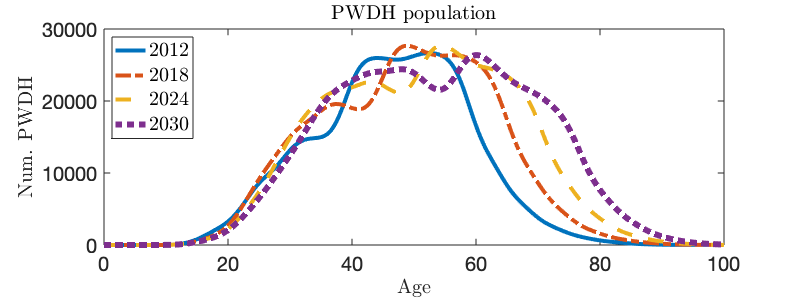}
    \caption{Evolution of the PWDH population age structure over the years 2012-2030. By year-end 2030, we see a bimodal distribution forming, with more PWDH around age 40 and age 60, as compared to PWDH around age 50. }
    \label{fig:PWHPrev2}
\end{figure}

 We remark, however, that not all age groups are projected to increase uniformly. As discussed in the previous subsections, future PWDH mortality rates are projected to decrease, however, such decreases are highly nonuniform based on age. Conversely, while new HIV diagnoses are expected to decrease, the projected age distribution does not show significant change. By 2030, this results in the formation of a bimodal-like population PWDH age distribtution, in which there are fewer PWDH approximately aged 50 compared to either age 40 and age 60 (Fig. \ref{fig:PWHPrev2}).

\section{Discussion}
In the present work, we have examined age-dependent data on PWDH in the United States to better understand how the developments over the previous 15 years have changed the age structure of the PWDH population in the United States and, furthermore, what implications these changes may have going forward. Starting with a simple age-structured population model, we used an Ensemble Kalman Filter Inversion, together with discrete age-structured HIV surveillance data, to reconstruct age-dependent PWDH mortality over the years 2008-22. We found that, due to the widespread availability of effective ART, mortality has declined significantly among PWDH since 2008, particularly among PWDH aged 40-80 years. While this trend reversed somewhat among PWDH aged 55 and older in 2021-22 due to the COVID-19 pandemic, data from 2022 suggest that COVID effects have subsided, likely signaling a return to pre-COVID trends.
\par We then developed a novel variant of Dynamic Mode Decomposition, nonnegative DMD, to develop projections of future changes in PWDH mortality and diagnoses. Projections of future PWDH mortality rates suggest that PWDH mortality will continue to decrease in the coming years, with the largest decreases expected among PWDH aged 40-55 and those older than 80. While mortality rates are expected to decrease among other age groups as well, smaller decreases are projected. We note that these projections assume a continuation of pre-existing trends, and cannot account for major changes in the coming years, such as medical breakthroughs.
\par Our projections of future HIV diagnoses suggest a slow, but notable, decrease in annual HIV diagnoses going forward, with around 28,000 diagnoses projected in 2030. Interestingly, our methods did not indicate a major change in the age structure of new HIV diagnoses, which has remained largely static in recent years, and is projected to remain so.
\par Applying our reconstructed and projected mortality and diagnosis data to our core PWDH age model, we then projected the age structure of the PWDH population over the coming years. Our simulations suggest that the PWDH population will age significantly. By 2030, the percentage of PWDH aged 55 years and older is expected to increase significantly (from 40\% to 47.4\%). Furthermore, these increases are more dramatic for older age groups; the percentage of PWDH aged 65 and older is expected to increase to 26.1\% by 2030, and the portion of those aged 75 years and older is expected to more than double by 2030.
\par We note that the analysis is subject to several important limitations. As previously mentioned, the projections of future changes in PWDH mortality and diagnoses assume continuation of pre-existing trends. It is possible that new medical technologies may alter the landscape of HIV prevention and care significantly. This is particularly important for new HIV diagnoses, as improvements such as rapid expansion of pre-exposure prophylaxis (PrEP) coverage or increased levels of viral suppression (for instance, due to long-acting injectable ART \cite{viguerie2025impact}), may result in fewer new HIV diagnoses than projected \cite{eisinger2018ending}. Furthermore, new HIV diagnoses were modeled as a linear source term, and any dependence on the PWDH age structure is assumed to be captured implicitly via any trends present in the data. However, the relationship between new the quantity and age distribution of new HIV diagnoses and the current PWDH population state may  be more complex, and future modeling efforts should seek to better explore and define this relationship.
\par Crucially, the current analysis represents only a first step towards developing comprehensive models for aging in the PWDH population. Future extensions should include other aspects of HIV prevention and care. For example, beyond biological age, time since acquiring HIV infection and time on ART may also be important factors in determining risk for certain commodities. Other extensions, including demographic stratification, may be similarly significant.
\par On the mathematical end, further analysis of nonnegative DMD is necessary. In the present work, we introduced the basics of the method and provided some elementary results for the current problem, as well as some intuitive arguments. Given the application-oriented nature of the current work, we elected to leave further exploration to future research. However, it is 
important to properly formalize how nonnegative DMD fits within the larger family of DMD methods. More rigorous analysis is necessary to establish when certain conditions observed in the current analysis, such as the existence of a limiting eigenvector with an associated simple eigenvalue of 1, can be expected to hold. Finally, developing more sophisticated and efficient computational approaches for nonnegative DMD, particularly those that exploit its parallelizability, are important for extending its application to larger-scale problems.
\par Our analysis highlights the urgent need to prepare for the aging population of PWDH. Current data suggests that PWDH are more vulnerable to aging-related comorbidities, requiring that the healthcare system be properly equipped to handle the potential rapid increase in such health conditions. Furthermore, in the case of preventable age-related comorbidities, the current analysis also emphasizes the need to develop and implement prevention programs to help reduce potential burden on the healthcare system as much as possible. Urgency is necessary, as model projections show that the portion of PWDH aged 75 and older will more than double over a five-year span. Beyond HIV care, using mathematical models to account for aging populations is an important direction for research in public health, and its importance is likely to increase in importance in the coming years.

\section*{Acknowledgments}
E. Iacomini is partially supported by the Italian Research Center on High-Performance Computing, Big Data and Quantum Computing (ICSC) funded by MUR Missione 4-Next Generation EU (NGEU) [Spoke 1 ``Future HPC \& Big Data"]. The authors are members of the INDAM GNCS (Italian National Group of Scientific Calculus. The authors would like to acknowledge Stefan Klus, Matthew Colbrook, Siobhan O`Connor and Ruiguang Song for their helpful input and advice.

\bibliographystyle{unsrt}  
\bibliography{references}  







\end{document}